\documentclass[nopreprintline,12pt,authoryear]{elsarticle}

\usepackage{amssymb}
\usepackage{amsmath}
\usepackage{subcaption}
\usepackage{makecell}

\setcellgapes{4pt}

\begin{document}

\begin{frontmatter}



\title{Retail Store Layout Optimization for Maximum Product Visibility}


\author[inst1]{Evren Gul\corref{corr}}
\author[inst1,inst2]{Alvin Lim}
\author[inst1]{Jiefeng Xu}

\address[inst1]{Nielsen Precima, LLC, 200 W.~Jackson Blvd., Chicago, IL 60606}
\address[inst2]{Emory University Goizueta Business School, 1300 Clifton Rd., Atlanta, GA 30322}
\cortext[corr]{corresponding author: evren.gul@nielseniq.com}

\begin{abstract}
It is well-established that increased product visibility to shoppers leads to higher sales for retailers. In this study, we propose an optimization methodology which assigns product categories and subcategories to store locations and sublocations to maximize the overall visibility of products to shoppers. The methodology is hierarchically developed to meet strategic and tactical layout planning needs of brick-and-mortar retailers. Layouts in both levels of planning are optimized considering eligibility requirements and complete set of shopper paths, thus, they successfully capture the unique shopping behavior of consumers in a store's region. The resulting mathematical optimization problem is recognized as a special instance of the well-known Quadratic Assignment Problem, which is considered computationally as one of the hardest optimization problems. We adopt a linearization technique and demonstrate via a real-world numerical example that our linearized optimization models substantially improve the store layout, hence, can be used in practical applications as a vital decision support model for store layout planning. 
\end{abstract}




\end{frontmatter}






\section{Introduction}

The retail industry is globally an essential sector which is known for slim profit margins and fierce competition. Grocery retail is additionally dynamic as it must answer to ever-changing customer shopping behaviors and business ecosystem. Therefore, grocers increasingly employ retail analytics to make strategic decisions in well-known processes such as pricing, assortment, marketing, and promotion planning. In addition to these traditional areas, store layout optimization can provide a retailer significant edge against the competition.
The present work focuses on store layout optimization to maximize visibility of products to shoppers as they navigate the store pick products. In this paper, we loosely associate product ``visibility'' and ``exposure'' and use the terms interchangeably.

The problem foundationally stems from the industry adage ``unseen is unsold''. Although grocery shoppers mostly rely on  predetermined shopping lists, unplanned purchases (also known as impulse purchases) can constitute 30-50\% of their purchases \citep{kollat1967customer, mishra2010impulsebuy}. Unplanned spending is triggered by in-store stimuli such as visual confrontations with a product \citep{rook1987buying}. Visual stimuli can be enhanced mainly by (1) increasing the shelf space allocated to products or (2) intelligently locating products in the store so that customers are exposed to maximum number of products along their shopping paths. \cite{dreze1994shelf} established that increased product visibility leads to higher impulse purchases. The first approach of visual stimuli was employed in \cite{xu2020maximizing} where floor space was optimally allocated to product categories such that the total store revenue is maximized. Related to the second approach, \cite{hui2013effect} assessed the hypothesis that traveling farther in a store results in higher unplanned spending. They statistically showed in a field experiment that as shoppers travel farther from their planned path, they make up to 50\% more unplanned purchases on average. In grocery industry, practitioners have long employed some tactics to lengthen shoppers’ travel in the store. A famous example of such tactic is ``hide milk at the back of the store''. These expertise-based strategies support layout decisions, yet, the final layout might be far from optimal because they do not determine the optimal layout based on how customers pick the products in their baskets. In this paper, we optimize the store layout via maximizing the total exposure of product categories and subcategories which are calculated over all customers' transactions.

In current practice, there are more design considerations to store layout than optimization to boost impulse buys. A store layout should provide plausible shopping experience through store atmosphere and in-store traffic patterns to enhance store loyalty and sales. Several researchers have shown the positive relationship between store layout and consumer purchase intentions \citep{baker1992experimental, merrilees2001superstore, ainsworth2017comfort}. Also, store layout should be aligned with in-store operations such as shelf stocking to provide operational efficiency \citep{lewison1994retailing, vrechopoulos2004virtual}. Consequently, well-designed store layouts appeal to the customers and contribute to both product sales and store profitability \citep{cil2012consumption}. In the present work, we aim to satisfy these desired layout characteristics and rules by introducing them as constraints in our mathematical formulation.

\begin{figure}[htb]
    \centering
    \includegraphics[width=5.9in,, trim = 0.8cm 16cm 0cm 2cm]
    {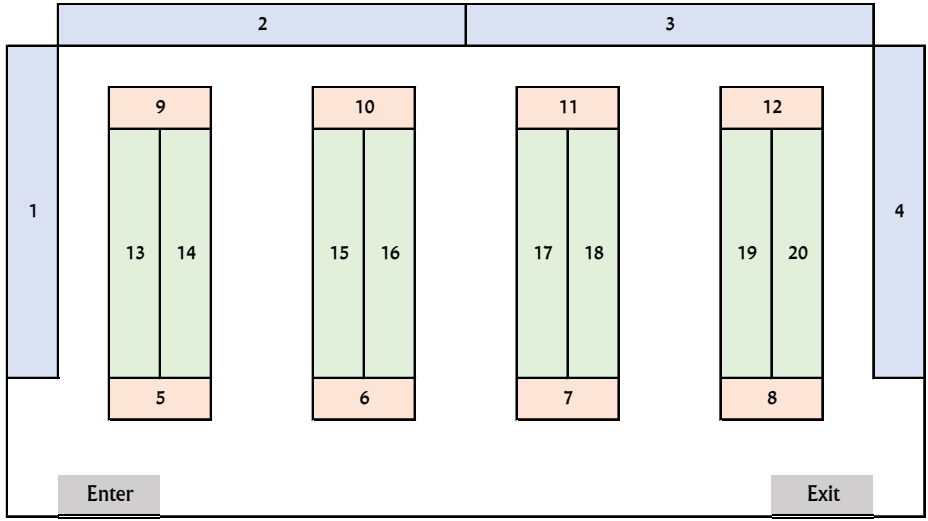}
    \caption{Grocery store with location level layout for category assignment.}
\label{fig:layout_lvl_1}
\end{figure}

\cite{vrechopoulos2004virtual} identified three major layouts in conventional retail settings: grid, free-form and racetrack. In grocery retail, we commonly encounter stores combining the grid and racetrack aisles as shown in Figure \ref{fig:layout_lvl_1}. The parallel arrangement of long aisles in the center of the store forms the ``grid'' which is circumvented by a ``racetrack'' with shelves leaning to the sides and back of the store. In this formation, main pathway on the perimeter branches into grid aisles, thus, it facilitates planned shopping behavior with flexibility and speed in finding products from a shopping list \citep{lewison1994retailing, levy2004retailing, larson2005exploratory}.

Retail facility layout optimization has only been in researchers’ radar in the recent decades. In the early work of \cite{botsali2005network}, different configurations of serpentine layout were analyzed to determine one which maximizes impulse purchase revenues. \cite{botsali2007retail} furthered the analysis by including the grid layout. \cite{cil2012consumption} built a customer-oriented supermarket layout using the product clusters derived from transactions. The layout is created to facilitate customer navigation by placing associated product clusters close to each other. The author claims that customer-appealing layouts would translate into higher purchase conversion rate.  Conversely, \cite{peng2011optimizing} disperses the defined ``must-have'' products on a grid layout in an aim to maximize the sales coming from ``impulse-buy'' products as customers travel between ``must-have'' items. \cite{yapicioglu2012retail} optimized a department store layout to maximize revenue with respect to the layout requirements. In the study, impulse purchases were reflected in revenues by scaling them with respect to the traffic densities of the store zones. The problem was solved with a Tabu Search algorithm on simulated problems of up to 20 departments on a racetrack structure. In their succeeding work, \cite{yapicioglu2012bi} translated the previous problem into a bi-objective model where they developed layout requirements as another objective of the problem. \cite{ozgormus2020data} also studied the similar bi-objective problem for a supermarket with racetrack aisles surrounding a grid. They simulated revenues on a function with predefined revenue coefficients and impulse purchase rates based on marketing and store manager's input. The problem is solved for 25 groups of product categories with a Tabu Search algorithm. In another implementation for a grocery store layout, \cite{flamand2016promoting} maximized the impulse revenues by determining the optimal locations of 31 product categories on different layout configurations. 

The aforementioned literature considers preset traffic densities as one of the drivers of the impulse purchases. However, traffic densities are calculated with respect to the locations of product categories in the current layout; these densities could potentially change with the optimized layout. Indeed, \cite{ballester2014effect} showed that traffic densities ultimately change with the changes in product category locations.  For instance, as a popular milk product is relocated in a supermarket, the associated traffic will follow it to the new location. To account for variable store traffic intrinsic to the problem, layout optimization can be carried out considering the shopper paths. As mentioned before, several researchers showed that product exposure, and eventually sales, increase with the increase in the length of shopping paths \citep{inman2009interplay, kholod2010influence, hui2013effect}. Supported by this finding, \cite{boros2016modeling} aimed to increase the path length for the average customer in a supermarket. They optimized the store layout for 23 departments and 27 representative customer baskets. Longer travel paths presumably increase the product exposure, however, product exposure can vary for different sections of the store. For instance, in Figure \ref{fig:layout_lvl_1}, product exposure is higher at the back and center of the store as products are displayed on both sides of the aisles, on the other hand, it is lower at the front of the store as products are only displayed on one side. Therefore, a better objective is direct optimization of product exposure. Recently, \cite{hirpara2021retail} maximized impulse purchase revenues due to departments' exposure along a customer's shopping path. They solved the problem for a store with 20 departments using a representative basket. 

The current problem, i.e., to find an optimal store layout to maximize exposure of product subcategories, is a real world application motivated by the store layout planning operations of a major grocery chain in Europe. It is studied on a typical medium-sized store located in an urban center. The representative layout, a combination of racetrack and grid structures, is provided in Figure \ref{fig:layout_lvl_1}. On the layout, 20 categories on three main shelf groups are defined with respect to fixture and product category requirements. These are peripheral shelves ($1-4$), endcaps ($5-12$) and aisle shelves ($13-20$). Based on this typical layout, store layout planners make decisions first at strategic and then at tactical levels. Strategic level store layout planning happens every five to ten years during which stores are reset for a new store concept and refreshed look. Stores are remodeled around established shopping behaviors and most recent industry trends. In our work, the strategic level layout model is termed the Level-1 model and it translates to finding 20 optimal category locations. In the model, categories are only allowed to move within the same prescribed locations due to fixture requirements. Optimal category locations from the Level-1 model form a basis for the tactical level layout planning which we address with a Level-2 model. 

\begin{figure}[htb]
    \centering
    \includegraphics[width=5.9in,, trim = 0.8cm 16cm 0cm 2cm]
    {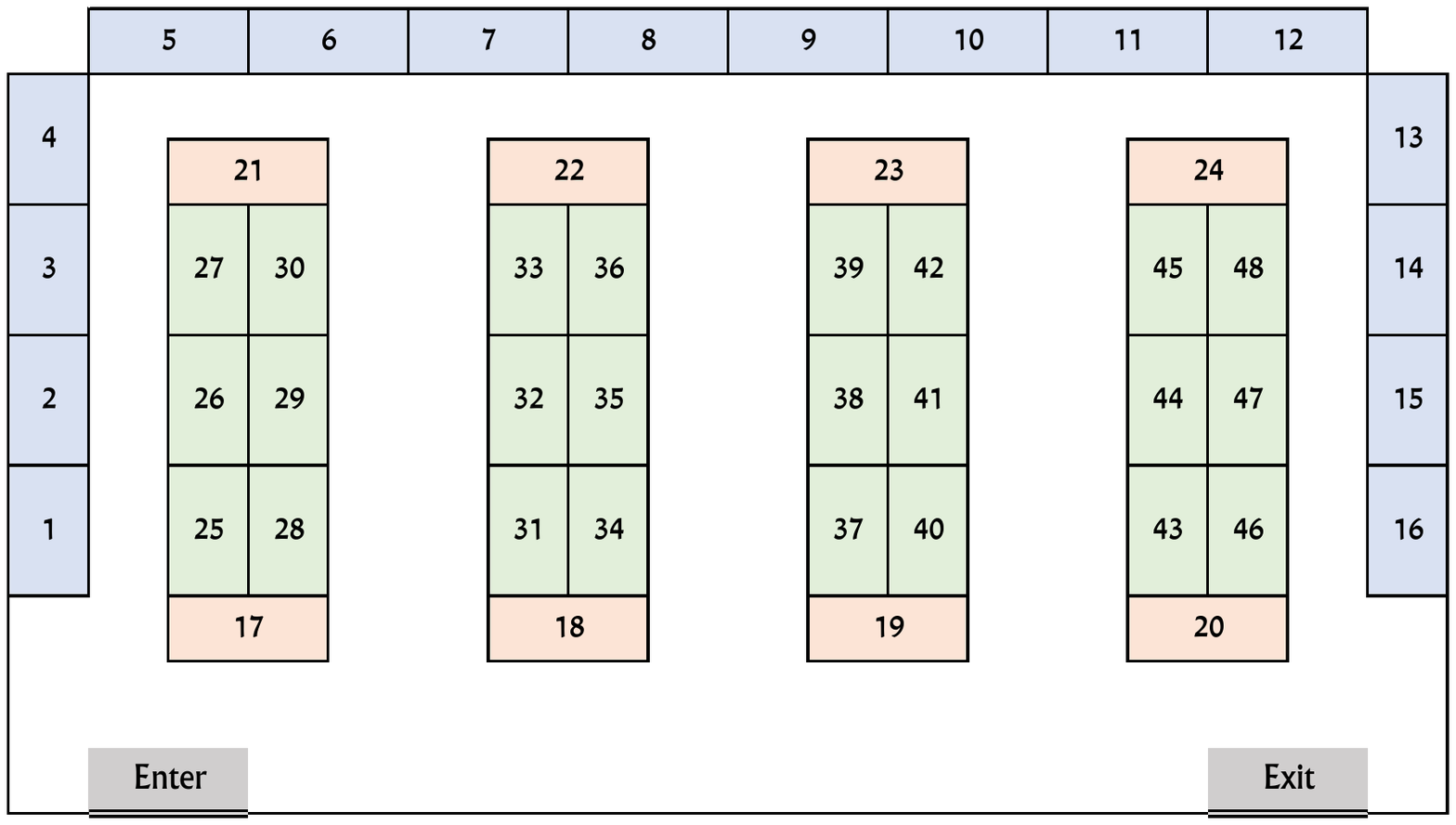}
    \caption{Grocery store with sublocation level layout for subcategory assignment.}
\label{fig:layout_lvl_2}
\end{figure}

The Level-2 model is used to assign subcategories on a more granular layout as in Figure~\ref{fig:layout_lvl_2} where a total of 48 subcategory sublocations are nested in the 20 category locations. The Level-2 problem determines the 48 optimal subcategory sublocations within already determined category locations from the Level-1 model. The optimization problem based on the Level-2 model is solved on a quarterly to yearly frequency so that store layout adapts to seasonal product assortment changes while maintaining the familiarity of the store's layout to the customers and keeping the operational costs to a minimum. Table~\ref{table:layout_summary} summarizes the assignment of product subcategories and their categories in the current store layout. In both the Level-1 and Level-2 problems, we aim to improve overall product exposures from the current store layout by optimally reassigning categories and subcategories to new locations to alter customers' shopping paths. 

The problems are solved for 20,842 customer transactions. We assume customers randomly pick the categories in their basket. Additionally, in the Level-2 problem, a customer is assumed to randomly pick all the subcategories belonging to a category before proceeding to the next category. Customers are considered to travel the shortest path between locations in the Level-1 problem (sublocations in the Level-2 problem). While finding the shortest path we consider center-to-center distance of the category locations and the subcategory sublocations in Level-1 and Level-2 problems, respectively. Next, we calculate the total exposure for each shortest path where exposure is defined as number of subcategory sublocations a customer passes by on her path.

The contribution of this paper is multi-fold. We first identify the shopper path-based maximum product exposure problems within the framework of the typical store layout planning processes for different planning horizons and formulate mathematical models for them. We further develop a mathematical programming-based solution approach to solve the computational intractable instances of the problems. We finally validate the algorithm using real customer transaction data and store layout location data and derive insightful outcomes. The rest of this paper is organized accordingly as follows. In Section~2, we develop mathematical optimization models for the store layout with maximum exposure problem required to support strategic level and tactical level planning as well as an integrated strategic-tactical planning model.  In Section~3, we present linearized versions of the proposed optimization models. Section~4 reports on our computational experience using a real-world planning instance from a European grocery chain.  We summarize our contribution and outline some future directions of research in Section~5.

\begin{table}[htbp]
\small
\caption{Product categories and their locations in a grocery store.} 
\centering 
\resizebox{\textwidth}{!}{\begin{tabular}{p{0.40\linewidth}  p{0.2\linewidth} p{0.48\linewidth}}
\hline 
Location - Category & Fixture/Shelf Type & Subcategory (Sublocation) \\ [1ex] 
\hline 
    1 - Fruits & Peripheral & Bananas (1), Citrus (2), Hard Fruit (3), Soft Fruit (4)  \\ 
    2 - Produce  & Peripheral & Root Vegetables (5), Leafy Green Vegetables (6), Ready Cut Vegetables (7), Fresh Juice (8) \\
    3 - Frozen Food & Peripheral & Frozen Vegetables (9), Frozen Snacks (10), Frozen Pizzas (11), Frozen Potatoes (12) \\
    4- Bakery & Peripheral & Bread Substitutes \& Spreads (13), Fresh Bread (14), Prepackaged Bread (15), Special Breads (16) \\
    5 - Fresh Dairy Food & Endcap & Fresh Dairy Food (17) \\
    6 - Seasonal Deli & Endcap & Seasonal Deli (18) \\
    7 - Seasonal Non-food & Endcap & Seasonal Non-food (19) \\
    8 - Seasonal Candy & Endcap & Seasonal Candy (20) \\
    9 - Fresh Dairy Drink & Endcap & Fresh Dairy Drink (21) \\
    10 - Fresh Nuts & Endcap & Fresh Nuts (22) \\
    11 - Household Assortment & Endcap & Household Assortment (23) \\
    12 - Seasonal Confectionery & Endcap & Seasonal Confectionery (24) \\
    13 - Appetizers & Aisle & Mediterranean Food (25), International (26), Seafood Appetizers (27) \\
    14 - Condiments & Aisle & Salad Essentials (28), Sauces (29), Pickled Food (30) \\
    15 - Preserved Food & Aisle & Preserved Vegetables (31), Preserved Meats (32), Sausages (33) \\
    16 - Chips & Aisle & Regular Chips (34), Premium Chips (35), Pop Corn (36) \\
    17 - Beer & Aisle & Bulk Beer (37), Bottle \& Can Beer (38), Special Beer (39) \\
    18 - Water \& Energy Drinks & Aisle & Energy Drinks (40), Flavored Water (41), Water (42) \\
    19 - Soda \& Soft Drinks & Aisle & Coke (43), Soda (44), Ice Tea (45) \\
    20 - Soups \& Eggs and Non-perishable Food & Aisle & Cookies (46), Soups (47), Eggs \& Non-perishable Dairy (48) \\ [1ex]
\hline 
\end{tabular}}
\label{table:layout_summary} 
\end{table}

\section{Model Formulation}
In this section, we first present an integrated model that covers both the Level-1 (strategic) and Level-2 (tactical) planning needs simultaneously. Due to practical reasons as well as the computational intractability of the integrated model, we then reformulate the integrated model into two models: Level-1 and Level-2. \\

First, we define some notations for index sets, constants and decision variables. \\

\hspace{-0.1in}\textbf{Index Sets:}
\begin{itemize}
  \item $K$: the set of all locations with the first and the last locations representing the entrance and the exit of the store, respectively.
  \item $K^*$: the set of all sublocations with the first and the last sublocations representing the entrance and the exit of the store, respectively.
  \item $I$: set of product categories with the first and the last categories representing customer check-in (a dummy ``product'' representing a customer going through the entrance location) and customer check-out (a dummy ``product'' representing traversal through the exit location), respectively.
  \item $I^*$: set of product subcategories with the first and the last subcategories representing customer check-in and check-out, respectively.
  \item $G_i$: set of product subcategories that belong to category $i$.
  \item $L_k$: set of sublocations that belong to location $k$.
\end{itemize}

\hspace{-0.1in}\textbf{Constants:}
\begin{itemize}
  \item $E_{k_1k_2}$: the exposure (number of sublocations passed by on the shortest path) from sublocation $k_1$ to sublocation $k_2$.
  \item $P_{i_1i_2}$: the number of occurrences (over all customers) for which a customer walks to pick
a product from subcategory $i_2$ immediately after picking a product from subcategory $i_1$. Dummy ``products'' such as customer check-in and customer check-out are treated like any other products in this context.
  \item $A_{ik}$: a binary constant with 1 indicating that product category $i$ is eligible to be assigned location $k$. For example, customer check-in is the only product eligible to be assigned the first location (the entrance of the store), and the frozen food category must be assigned only to locations where freezers can be located.
\end{itemize}

\hspace{-0.1in}\textbf{Decision Variables:}
\begin{itemize}
  \item $x_{ik}$: the binary variable with value $1$ if and only if product category $i$ is assigned to store location $k$, and $0$ otherwise.
  \item $z_{i_1k_1}$: the binary variable with value $1$ if and only if product subcategory $i_1$ is assigned to sublocation $k_1$, and $0$ otherwise.
\end{itemize}

With the above notations, we can now formulate the Integrated Model that covers both category assignments at the strategic level and subcategory assignments at the tactical level simultaneously.\\

\hspace{-0.1in}{\bf Integrated Model (IM):}
\begin{alignat}{3}
  \textrm{Maximize} \quad &\sum_{k_1,k_2 \in K^*} \sum_{i_1,i_2 \in I^*} E_{k_1k_2}P_{i_1i_2} z_{i_1k_1} z_{i_2k_2} && \label{eqn:obj1}\\
 \textrm{subject to}\quad & \sum_{i \in I} x_{ik} =1, &\quad &\forall k \in K, \label{eqn:assign_loc1}\\
  &\sum_{k \in K}  x_{ik} = 1,  &\quad & \forall i \in I, \label{eqn:assign_cat1}\\
 &x_{ik} \leq A_{ik},  &\quad & \forall i \in I, k \in K, \label{eqn:assign_eligible1}\\
 &\sum_{k_1 \in L_k} z_{i_1k_1} = x_{ik}, &\quad & \forall i \in I, k \in K, i_1 \in G_i, \label{eqn:subassign_loc1}\\
&\sum_{i_1 \in G_i} z_{i_1k_1} = x_{ik}, &\quad & \forall i \in I, k \in K, k_1 \in L_k, \label{eqn:subassign_cat1}\\
  &x_{ik}  \in \{0,1\}, &\quad & \forall i \in I, k \in K, \label{eqn:binary_var_x1}\\
  &z_{i_1k_1}  \in \{0,1\}, &\quad & \forall i_1 \in I^*, k_1 \in K^*. \label{eqn:binary_var_z1}
\end{alignat}

In the above formulation, the objective function \eqref{eqn:obj1} maximizes the total exposure across all customers' shopping paths. Constraints \eqref{eqn:assign_loc1} and \eqref{eqn:assign_cat1} stipulate the one-to-one assignment requirement of a category to a location, more specifically, only one category can be assigned to a specific location, and only one location can be assigned to a specific category. Constraint \eqref{eqn:assign_eligible1} specifies that the assignment of a category to a location must be feasible subject to practical eligibility conditions. Constraint \eqref{eqn:subassign_loc1} and \eqref{eqn:subassign_cat1} establish further requirements in a Level-2 (tactical level) layout that a product subcategory must be assigned to a sublocation within the location that its corresponding category is assigned to. Finally, constraints \eqref{eqn:binary_var_x1} and \eqref{eqn:binary_var_z1} guarantee the binary nature of the decision variables.

The IM can be considered as a special instance of the {\em Quadratic Assignment Problem} (QAP) \citep{koopmans1957assignment}.  The objective of the standard QAP is to assign a set of facilities to a set of locations in such a way that the total assignment cost is minimized. In our problem, we can view each product subcategory as a facility and construct a maximization instance of the QAP with additional restrictions imposed by the eligibility constraints \eqref{eqn:assign_eligible1},  \eqref{eqn:subassign_loc1} and \eqref{eqn:subassign_cat1}. QAP is a well-known {\em NP-hard} combinatorial optimization problem \citep{sahni1976p}. Although extensive research has been done for over the last six decades, the QAP is still recognized computationally as one of the hardest known optimization problems. No known exact algorithm can solve practical problems of sizes more than 30 locations in reasonable computational time. For example, \cite{brixius2000solving} reported an interesting computational experience in solving a 30-location benchmark QAP (Nugent 30 from QAPLIB) using a branch-and-bound method on grid computers. It required 1000 workstations over seven days to solve the QAP instance. Our preliminary computational experiments also confirmed that the integrated model with 48 sublocations/subcategories cannot be solved directly using a state-of-the-art quadratic programming solver such as Gurobi.

From a practical implementation perspective, as mentioned earlier, the strategic level store layout problem (with the Integrated Model as an instance) is solved infrequently (every five to ten years) and is only required for a store reset. On the other hand, the tactical level problem is required to be solved as often as quarterly to adjust the subcategory-sublocation assignments based on stationery category location assignments to accommodate seasonal adjustments to product assortment. This means that, from a practical standpoint, it is much more important for instances of a tactical level problem to be solved efficiently than a strategic level problem such as the Integrated Model. For this reason, we reformulate the Integrated Model into two distinct models:
\begin{itemize}
    \item {\bf Level-1 Model:} Optimally assign categories to locations; and 
    \item {\bf Level-2 Model:} Optimally assign subcategories to sublocations given category locations from a Level-1 Model.
\end{itemize}

We define some additional constants associated with the Level-1 Model:
\begin{itemize}
\item $\bar{E}_{k_1k_2}$: the exposure (number of sublocations passed by on the shortest path) from location $k_1$ to location $k_2$. At this level, we assume all the sublocations of category is exposed when category is picked.
  \item $\bar{P}_{i_1i_2}$: the number of occurrences (over all customers) for which a customer walks to pick a product from category $i_2$ immediately after picking a product from category $i_1$.
\end{itemize}

\noindent Then the Level-1 Model can be formulated as follows:\\

\hspace{-0.1in}{\bf Level-1 Model (L1M):}
\begin{alignat}{3}
  \textrm{Maximize} \quad &\sum_{k_1,k_2 \in K} \sum_{i_1,i_2 \in I} \bar{E}_{k_1k_2}\bar{P}_{i_1i_2} x_{i_1k_1} x_{i_2k_2} && \nonumber \\
 \textrm{subject to} \quad & \sum_{i \in I} x_{ik} =1, &\quad &\forall k \in K, \label{eqn:assign_loc2}\\
  &\sum_{k \in K}  x_{ik} = 1,  &\quad & \forall i \in I, \label{eqn:assign_cat2}\\
 &x_{ik} \leq A_{ik},  &\quad & \forall i \in I, \quad k \in K, \label{eqn:assign_eligible2}\\
   &x_{ik}  \in \{0,1\}, &\quad & \forall i \in I, \quad k \in K. \label{eqn:binary_var_x2}
\end{alignat}

Letting $\{\bar{x}_{ik}\}_{i \in I, k \in K}$ be an optimal solution obtained from L1M, we formulate the corresponding Level-2 Model as follows:\\

\hspace{-0.1in}{\bf Level 2 Model (L2M):}
\begin{alignat}{3}
  \textrm{Maximize} \quad &\sum_{k_1,k_2 \in K^*} \sum_{i_1,i_2 \in I^*} E_{k_1k_2}P_{i_1i_2} z_{i_1k_1} z_{i_2k_2} && \nonumber \\
 \textrm{subject to}\quad & \sum_{k_1 \in L_k} z_{i_1k_1} = \bar{x}_{ik}, &\quad & \forall i \in I, k \in K, i_1 \in G_i, \label{eqn:subassign_loc3}\\
&\sum_{i_1 \in G_i} z_{i_1k_1} = \bar{x}_{ik}, &\quad & \forall i \in I, k \in K, k_1 \in L_k, \label{eqn:subassign_cat3}\\
  &z_{i_1k_1}  \in \{0,1\}, &\quad & \forall i_1 \in I^*, \quad k_1 \in K^*. \label{eqn:binary_var_z2}
\end{alignat}
\noindent

We expect L1M to have higher objective function value than L2M as exposure values ($\bar{E}_{k_1k_2}$) in L1M are greater than ($E_{k_1k_2}$) in L2M. This is due to the fact that we assume all subcategories of a category is exposed when passed by in L1M, on the other hand, only the subcategories passed by are counted in L2M. Therefore, L1M broadly maximizes total exposure by optimally determining category locations with allowed category movements on the whole store footprint, and L2M improves total exposure by locally assigning the sublocations to subcategories with movements restricted to their own category location.
Additionally, note that in the model L2M, $\bar{x}_{ik}$ for all $i \in I, k \in K$ are constants rather than decision variables. Also, just like IM, both L1M and L2M are instances of the QAP with restrictions and are still computationally challenging to solve.  In the next section, we discuss a linearization technique that will allow us to significantly reduce the computing time. 

\section{Linearization of the QAP}
As described above, the store layout with maximum exposure problem can be formulated as an instance of the QAP and modeled as a standard mixed integer quadratic program (MIQP). It is well-known that the quadratic terms in the objective function can be eliminated and transformed into linear terms by introducing new variables and new linear (and binary) constraints \citep{watters1967reduction}. After the linearization, the MIQP becomes a mixed (0-1) integer linear programming problem (MILP), so the existing methods and standard MILP solvers can be applied to solve the problem. Although, in general, MILP is considered as an ``easier'' problem than its counterpart MIQP, but the very large number of new variables and constraints increase the original problem size significantly leading to potentially considerable requirement in additional computing resources.

Researchers developed a variety of linearization techniques to transform an MIQP into an MILP depending various methods to replace the quadratic terms. Some of the better known linearization methods include those proposed by \cite{lawler1963quadratic}, \cite{glover1974converting}, \cite{kaufman1978algorithm}, \cite{frieze1983quadratic}, and \cite{adams1994improved}. After experimenting with these linearization techniques, we adopt the Adams-Johnson linearization method from \citep{adams1994improved} for the much reduced solution times of the resulting MILP formulation as compared to the original MIQP model. 

For illustrative purposes, we show how we apply the Adams-Johnson relaxation method to L1M. The same procedure can be applied to IM and L2M directly to obtain associated MILPs. For each quadratic product term $x_{i_1k_1}x_{i_2k_2}$ in the objective function of L1M, we introduce a continuous variable $w_{i_1k_1i_2k_2}$ with the following additional constraints:
\begin{alignat}{2}
\sum_{i_1 \in I} w_{i_1k_1i_2k_2} =x_{i_2k_2}, \quad  \quad&\forall i_2 \in I, \quad  k_1,k_2 \in K, \label{eqn:linear1}\\
 \sum_{k_1 \in K}  w_{i_1k_1i_2k_2} = x_{i_2k_2},  \quad  \quad& \forall i_1, i_2 \in I, \quad  k_2 \in K, \label{eqn:linear2}\\
  w_{i_1k_1i_2k_2}=w_{i_2k_2i_1k_1}, \quad  \quad& \forall i_1, i_2 \in I, \quad  k_1,k_2 \in K, \label{eqn:linear3}\\
  w_{i_1k_1i_2k_2}  \geq 0, \quad \quad & \forall i_1, i_2 \in I, \quad k_1,k_2 \in K. \label{eqn:y_var}
\end{alignat}

\noindent Thus, the linearized L1M becomes\\

\hspace{-0.1in}{\bf Linearized Level-1 Model (LL1M):}
\begin{alignat}{2}
\noindent
   \textrm{Maximize} \quad \quad &\sum_{k_1,k_2 \in K} \sum_{i_1,i_2 \in I} \bar{E}_{k_1k_2}\bar{P}_{i_1i_2} w_{i_1k_1i_2k_2}  \nonumber \\
 \textrm{subject to}\quad \quad & \eqref{eqn:assign_loc2},\eqref{eqn:assign_cat2}, \eqref{eqn:assign_eligible2},\eqref{eqn:binary_var_x2},\eqref{eqn:linear1},\eqref{eqn:linear2},\eqref{eqn:linear3},\eqref{eqn:y_var}, \nonumber
\end{alignat}

\noindent and for completeness, we present the linearized versions of the L2M and IM as follows.\\

\hspace{-0.1in}{\bf Linearized Level-2 Model (LL2M):}
\begin{alignat}{3}
   \textrm{Maximize} \quad &\sum_{k_1,k_2 \in K} \sum_{i_1,i_2 \in I} E_{k_1k_2}P_{i_1i_2} y_{i_1k_1i_2k_2} \nonumber  &\\*
 \textrm{subject to} \quad & \eqref{eqn:subassign_loc3},\eqref{eqn:subassign_cat3},\eqref{eqn:binary_var_z2}, &\nonumber\\*
&\sum_{i_1 \in I} y_{i_1k_1i_2k_2} = z_{i_2k_2},&\quad & \forall i_2 \in I^*, k_1,k_2 \in K^*, \label{eqn:linear1_2}\\*
& \sum_{k_1 \in K}  y_{i_1k_1i_2k_2} = z_{i_2k_2}, &\quad  & \forall i_1, i_2 \in I^*, k_2 \in K^*, \label{eqn:linear2_2}\\*
 & y_{i_1k_1i_2k_2}=y_{i_2k_2i_1k_1}, &\quad & \forall i_1, i_2 \in I^*, k_1,k_2 \in K^*, \label{eqn:linear3_2}\\
 & y_{i_1k_1i_2k_2}  \geq 0, &\quad & \forall i_1, i_2 \in I^*, k_1,k_2 \in K^*. \label{eqn:y_var_2}
\end{alignat}
  
\hspace{-0.1in}{\bf Linearized Integrated Model (LIM):}
\begin{alignat}{2}
\noindent
   \textrm{Maximize} \quad \quad &\sum_{k_1,k_2 \in K} \sum_{i_1,i_2 \in I} E_{k_1k_2}P_{i_1i_2} y_{i_1k_1i_2k_2}  \nonumber \\
 \textrm{subject to}\quad \quad & \eqref{eqn:assign_loc1}, \eqref{eqn:assign_cat1},\eqref{eqn:assign_eligible1},\eqref{eqn:subassign_loc1}, \eqref{eqn:subassign_cat1},\eqref{eqn:binary_var_x1}, \eqref{eqn:binary_var_z1},   \eqref{eqn:linear1_2},\eqref{eqn:linear2_2},\eqref{eqn:linear3_2}, \eqref{eqn:y_var_2}. \nonumber
\end{alignat}

\section{Computational Results}
In this section, we report our experience solving a real instance of a store layout with maximum exposure problem for a major grocery chain in Europe. Although, in practice, we only need to solve a Level-2 problem in most situations as category locations were fixed, we report our experience solving each of IM, L1M and L2M for completeness.

We conduct the computational experiments based on real data for a typical medium-sized store of the grocery chain mentioned in Section~1. As shown in Figure 1 and 2, the problem involves 20 locations/categories (22 nodes including the store entrance and the store exit) at Level-1 and 48 sublocations/subcategories (50 nodes including the store entrance and store exit) at Level-2, respectively. We employ Gurobi (version 7.02), an industry-leading mathematical programming solver package to solve our store layout problems directly using its MIQP or MILP solvers. We report results for solving models specified in IM, L1M and L2M, as well as their linearized counterparts LIM, LL1M and LL2M. The computational tests were conducted using an Amazon Elastic Compute Cloud (EC2) server (running Amazon Linux AMI 2018.03 with dual Intel\textsuperscript{\textregistered} Xeon\textsuperscript{\textregistered} 16-core E5-2686 v4 CPU with 2.3GHz and 244GB memory). The number of computing threads for Gurobi is set to 16. All other Gurobi parameters are set to their default values except for those described thereafter. 

We first report our experiences on solving the Integrated Model directly as shown in Table~\ref{table:im} below.  As expected, directly solving the integrated model resulted in intractable outcomes consistent with its QAP counterparts as reported in the literature. After 20,160 minutes (two weeks) of computation time, Gurobi yields a feasible solution for the IM with a 1,609\% optimality gap and an objective function value (OFV) of 346,165. For the linearized model LIM, Gurobi even failed to solve the linear programming relaxation at the root node of the branch and bound tree after 20,160 minutes of computing time.

\begin{table}[htbp]
\centering \makegapedcells
\caption{Computational Results on the Integrated Models IM and LIM.}
\label{table:im}
\begin{tabular}{|l|c|c|}
\hline
Result \textbackslash Model& IM& LIM\\
\hline
Solution OFV &346,156 &No feasible solution\\
\hline 
Optimality Gap (\%) &1,609 &N/A\\
\hline
CPU Time (minutes) & 20,160 &20,160\\
\hline
\end{tabular}
\end{table}

Now, we report some results using a two-level hierarchical approach, i.e., we solve the strategic layout decision problem L1M using aggregate data first, then based on the optimal category/location layout obtained from L1M, we solve the tactical subcategory/sublocation layout problem with L2M subsequently. It is worth noting that the Level-2 models L2M or LL2M can start with any feasible solution at Level-1 (e.g., the current layout at the strategic level) directly, and are applied as stand-alone optimizations in most cases in real-world applications. To avoid the numerical issues that may obscure the performance of Gurobi, we set both the tolerances for optimality gap and feasibility gap to $1e{-}4$.  We summarize in Table~\ref{table:two_level} the results of solving L1M and L2M, as well as their linearized counterpart models, LL1M and LL2M.

\begin{table}[htbp]
\centering \makegapedcells
\caption{Computational Results on the Two Level Models.}
\label{table:two_level}
\begin{tabular}{|l|c|c|c|c|}
\hline
Result \textbackslash Model& L1M& L2M&LL1M&LL2M\\
\hline
Solution OFV &1,177,636 &811,540&1,177,636&806,342\\
\hline
Optimality Gap (\%)&0.0067&0.0099 &0.0000&0.0000\\
\hline
CPU Time (minutes) &15.2 &516.5 &0.1&2.4\\
\hline
\end{tabular}
\end{table}

We observe from Table~\ref{table:two_level} that Gurobi is able to solve the L1M and L2M problems optimally using the hierarchical sequential approach. It takes 15.2 minutes to solve the L1M and 516.5 minutes to solve the subsequent L2M. The OFV from the final optimal solution of L2M is 811,540, which is far superior to the feasible solution (with OFV 346,156) obtained by Gurobi from the IM. Table~\ref{table:two_level} also clearly shows the computational time improvement from the linearized models LL1M and LL2M. Gurobi takes just 0.1 minutes to solve LL1M and 2.4 minutes to solve LL2M, respectively. This indicates that the linearized versions of the layout decision models, regardless of tactical or strategic levels, can dramatically reduce the computation time by orders of magnitude. This makes the approach of using  LL1M and LL2M attractive as a vital tool in making store layout decisions in practice.

Moreover, we observe that the final OFVs from L2M and LL2M are different: 811,540 versus 806,342, with the same OFV at 1,177,636 from the associated Level-1 models L1M and LL1M. However, their optimality gaps from L2M and LL2M indicate that both solutions are optimal or very close to the optimal. This is because in this particular case, L1M and LL1M yield two alternative optimal solutions with the same OFV. Sequentially, L2M and LL2M built upon two different models based on these two Level-1 solutions, and yield two ultimate different optimal solutions at the tactical level. 

Although we emphasize that L2M or LL2M are the models required to make short-term layout decisions needed in most applications in practice, there are a few circumstances where combined strategic and tactical layout decisions need to be made simultaneously (e.g., in a store remodeling). Since it is impractical to employ IM or LIM, we recommend the adoption of an approximate method like the two level hierarchical sequential approach described above. As discussed earlier, such an approach of decomposing the original problem into two sub-problems then solving them sequentially may lead to a sub-optimal final solution for the original problem. To alleviate this shortcoming, we devise a naive heuristic that simply allows Gurobi to generate a pool of optimal or near optimal solutions for L1M/LL1M. Then for each solution in the resulting solution pool, we generate the corresponding L2M or LL2M model and solve it to optimality. Finally, we select the best solution from these pseudo-optimal solutions as the final best solution. We implement such a mechanism in our two-level hierarchical framework by allowing 10 optimal/near optimal solutions in the solution pool for L1M/LL1M. Note that the solutions in the pool are selected so that the OFVs are all within 0.1\% of the best solution. The results are summarized in Table~\ref{table:pooled_method}.

\begin{table}
\centering \makegapedcells
\caption{Computational Results for Pooled Heuristic Using Pool Size 10.}
\label{table:pooled_method}
\begin{tabular}{|l|c|c|}
\hline
Results \textbackslash Model& L1M$\rightarrow$L2M& LL1M$\rightarrow$LL2M\\
\hline 
Level-1 Worst/Best OFV & 1,177,613/1,177,636 & 1,177,520/1,177,636\\
\hline
Level-2 Worst/Best OFV & 796,562/813,664 & 795,677/813,664\\
\hline
Final Best OFV &813,664 & 813,664\\
\hline
CPU Time (minutes) & 11,922 &42\\
\hline
\end{tabular}
\end{table}

Table~\ref{table:pooled_method} confirms that the proposed pooled method is effective in improving the final sub-optimal solutions for both L1M$\rightarrow$L2M and LL1M$\rightarrow$LL2M sequential approaches. Compared to the results from Table~\ref{table:two_level} (which is the case with the pool size 1), using a pool size of 10 improves the final OFV from 811,540 to 813,664 for L1M$\rightarrow$L2M and from 806,342 to 813,664 for LL1M$\rightarrow$LL2M. The total computation time using LL1M$\rightarrow$LL2M is 42 minutes, which is dramatically less than the 11,922 minutes used for L1M$\rightarrow$L2M. Particularly noteworthy is that the further improvement in the quality of the solution from LL1M$\rightarrow$LL2M using the pooled heuristic was obtained with very affordable incremental computational time. 

\begin{figure}[htb!]
\captionsetup[subfigure]{justification=centering}
\centering
\begin{subfigure}[b]{0.55\textwidth}
\centering
   \includegraphics[width=4.75in,, trim = 5.5cm 16cm 0cm 2cm]{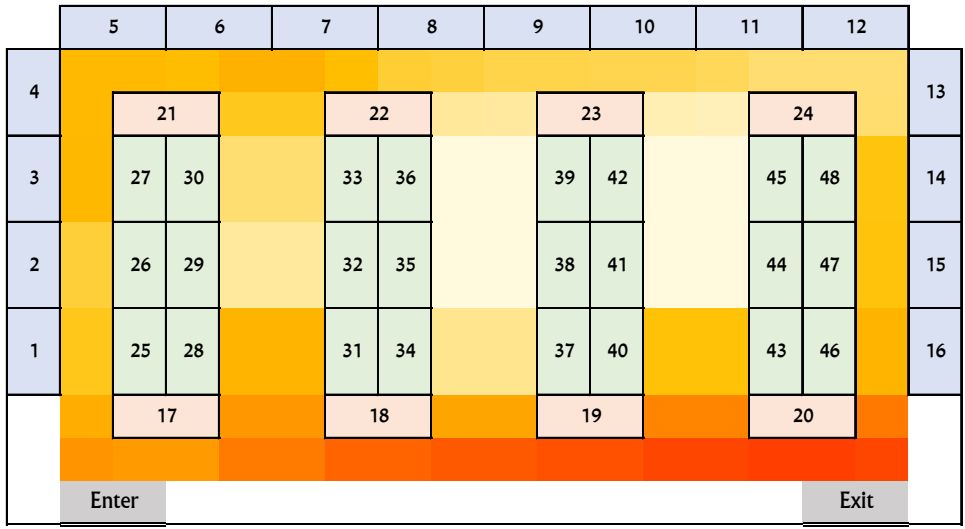}
   \vspace{-2\baselineskip}
   \caption{Current layout.}
   \label{fig:current_layout} 

\end{subfigure}

\begin{subfigure}[b]{0.55\textwidth}
\centering
   \includegraphics[width=4.75in,, trim = 5.5cm 16cm 0cm 0.75cm]{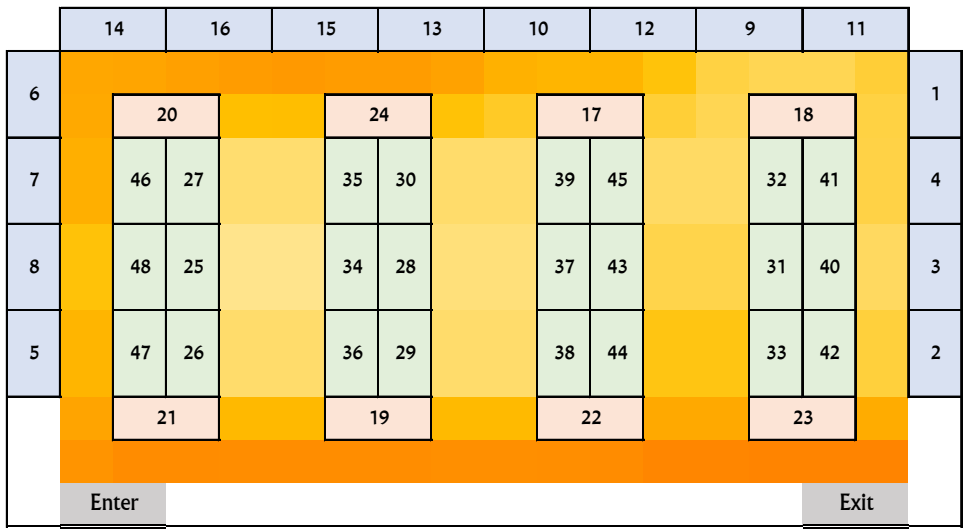}
   \vspace{-2\baselineskip}
   \caption{Optimal layout.}
   \label{fig:optimal_layout}
\hspace{5pt}
\end{subfigure}

\caption{Subcategory assignments and traffic density maps for (a) Current and (b) Optimal layouts.}
\label{fig:traffic_dens}
\end{figure}

We provide the current and optimal layouts with subcategory assignments and associated traffic density maps in Figure~\ref{fig:traffic_dens} on which traffic density increases from pale yellow to red. The optimal layout is the output of the two-level hierarchical approach discussed above. We achieve significant improvement in product exposure in the optimal layout where total product exposure is increased by 9.4\% as compared to the current layout. Yet, the total traveled distance increased only by 5.4\%. Field studies by \cite{west1951impulserate1} and  \cite{bellenger1978impulserate2} on impulse purchase percentages for different product groups serve as support points for rate generation in the literature. In these field studies, impulse purchase rates range from 27\% to 69\%, therefore, it is expected that 9.4\% more visibility of products to shoppers in the optimal layout will result in substantial incremental revenue from unplanned purchases.

Figure~\ref{fig:current_layout} shows that traffic on the current layout is highly concentrated at the front of the store where the exposure is lower as subcategories are only located on one side. On the other hand, in the optimal layout in Figure~\ref{fig:optimal_layout}, traffic density is much more distributed as compared to the current layout. Some eminent observations may provide explanations on how total exposure is improved. First, the ``Breads'' category (13-16), one of the top selling categories, is moved from right-hand side of the store to the back of the store in the Level-1 optimization. As a result, it leads to more traffic and exposure deeper into the store as customers have to pass by the grid aisles to pickup bread products. Secondly, ``Bananas'' (1), the most selling subcategory within ``Fruits'' category (1-4), is located close to the entrance at the bottom left corner of the store in the current layout. In Level-1 optimization, first, the ``Fruits'' category is moved to the peripheral shelves at the right-hand side of the store. Then, in the Level-2 optimization, ``Bananas'' moved to the far corner of the store which in return further increases exposure. Finally, 7 out of 8 subcategories (25, 28, 31, 34, 37, 40, 43, 46) , which are located in the section of grid aisles that is nearer the front of the store in the current layout, are moved to the center of grid aisles in Level-2 optimization. Consequently, traffic concentration at the front of the store in the current layout is dispersed to the store center along with more exposure in the optimal layout. Overall, the Level-1 problem determines the optimal category locations which maximize the total exposure and Level-2 optimization further enhances the total exposure by optimal arrangement of the subcategories within the categories. 




\section{Conclusion}
In this research, we propose a data-driven optimization methodology for retail store layout. 
This paper contributes to the literature with a novel optimization framework which can support layout planning operations over a wide time span along with maximization of total product visibility on the layout through well-depicted customer mobility using shopper paths. Specifically, this new approach optimizes the store layout while conforming to retailer's layout operations and decisions over entire planning horizon and dynamically accounting for changes in product exposures due to resulting changes in store traffic with the layout change. Hierarchical optimization models, Level-1 at the strategic level and Level-2 at the tactical level, integrate into decision processes for both long and short terms. Optimization models are not only in-line with decision processes but also satisfy additional practical business constraints such as suitability of fixtures. In both models, category and subcategory assignments to optimal locations maximize overall product exposure by inducing customer shopping paths that maximize traversal of the store layout. Optimization over customer shopping paths provides an edge to the approach as the store automatically adapts its layout to the changes in the customer habits and industry trends. 

We formulate the integrated store layout problem at both strategic and tactical levels as a restricted case of the hard-to-solve quadratic assignment problem. Furthermore, we decompose the integrated problem into the two sub-problems at the strategic and tactical levels, respectively. The Adams-Johnson linearization technique is used to expedite the solution process by orders of magnitude. We show that the real world problems at the strategic level as well as the tactical level can be readily solved using a standard solver package like Gurobi. In addition, the pooled heuristic coupled with the linearization technique is proven effective for solving the integrated model when needed. Results using the data set from a real grocery store clearly show the benefit in using the approach as total exposure is significantly increased while satisfying operational requirements. 

The research described in this paper reflects our initial effort to address the retail store layout decisions based on data-driven analytics and optimization methods. The problem set up in this paper provides a baseline for retail store layout optimization on which several extensions of the problem can consider. First, we assume categories and subcategories can be interchanged because they have equal size. As an extension, problem can be solved by addressing complications of size and interchangeability constraints and different store layouts. Secondly, the problem can be transformed into an impulse revenue maximization problem by converting product exposures to revenues using impulse purchase rates or expert opinion as discussed in the literature review. Another natural extension would be to combine the problem with the shelf/floor space allocation problem. We believe the latter two extensions would be more meaningful with rigorous estimation of impulse revenues.

As a final note, the methodology proposed can be flexibly configured to serve different business objectives. Retailers desire similar layouts for stores in a region to retain familiarity and easy navigation as customers may visit different stores. In-line with this requirement, first, the strategic Level-1 problem can be solved for a group of stores with similar layouts for pooled shopper paths. The resulting optimal layout would both provide the uniformity and maximize the total exposure for all stores. Next, each store's layout can be customarily optimized considering their own shopper paths which further boost product exposure without harming the store identity. 

\bibliographystyle{elsarticle-harv}
\bibliography{arXiv_-_optimal_store_layout_for_maximum_exposure}

\end{document}